\documentclass[12pt]{article}

\usepackage[ansinew]{inputenc}
\usepackage{color}
\usepackage{enumerate,latexsym}
\usepackage{amsmath,amssymb}
\usepackage{graphicx}
\usepackage{amsmath,amssymb}
\usepackage{amsxtra}
\usepackage{amsfonts}
\usepackage{enumerate}
\usepackage{wrapfig}
\usepackage{xlop}
\usepackage{bbm}

\DeclareMathOperator*{\argmin}{arg\,min}

\newtheorem{theorem}{Theorem}[section]

\usepackage[colorlinks = true,
            linkcolor = blue,
            urlcolor  = blue,
            citecolor = blue,
            anchorcolor = blue]{hyperref}
\begin{document}

\title{Random Graph Models and Matchings}
\author{Lucas G. Rooney\\
College Natural Sciences,\\
University of Massachusetts,\\
Amherst,
MA,
01003\\
\texttt{lucasrooney@umass.edu}}
\date{May 9, 2019}
\maketitle

\begin{abstract}
In this paper we will provide an introductory understanding of random graph models, and matchings in the case of Erd\H{o}s-R\'enyi random graphs.  We will provide a synthesis of background theory to this end.  We will further examine pertinent recent results and provide a basis of further exploration.
\end{abstract}

\section{Introduction}
Random graph models are, in the most basic sense, means by which to construct random graphs which synthetically emulate the topology of real-world networks.\cite{Vince}  An intuitive appreciation of this directly follows---if we can build the topology of a real-world network into a tractable random graph model, then we can gain a richer and more accurate understanding of the characteristics of that network.  An important application of random graph modeling is to random graph matchings.

Graph matching is a rich area of statistics literature, particularly as the revolution in computing has been taking place in the last several decades.\cite{Vince}  Applications of graph matching include, but are not limited to, computer vision, pattern recognition, manifold and embedded graph alignment, shape matching and object recognition.\cite{Vince}  Graph matching is also a special case of the quadratic assignment problem, which is NP-hard, and given certain parameters being permitted in the generated random graphs, is actually equivalent to that problem.\cite{Vince}  While an intuitive understanding of graph matching as structure-preserving alignment between graphs is fairly accessible, we seek to build a rigorous backbone on which a technical understanding can be accomplished.

\section{Background}
\subsection{Networks}
In the intuitive sense, networks are the fabric which underlie systems in the world around us---our circles of friends, digital communications, political affiliations and economic processes all are examples of networks.  More technically, networks comprise neurological connections, molecular structures, and biochemical processes.  Mathematically, though, networks can be thought of as graphs.\cite{Eric}  Namely, a node in the network can be thought of as a vertex of a graph, and an edge connecting two vertices represents some communication between the two nodes.  We can see this with an example from Zachary's karate club, where nodes represent members of a karate club, and two nodes being linked means that the two members interacted.  This particular model is highly cited since it provides a strong basis for understanding subgroupings as, due to a fissure in the social structure in the club, there are two known subgroups within the 34 members.\footnote{There are other informative and interesting visualizations of networks, but given the scope of this paper they are not included.  Kolaczyk provides a comprehensive analysis of network data and structures for further reading.\cite{Eric}}

\begin{figure}
\centering
\includegraphics[width=.5\textwidth]{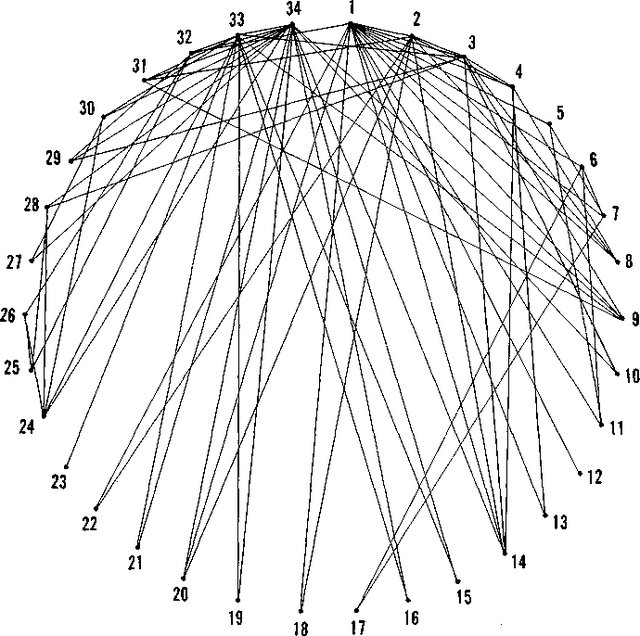}
\caption{This is a graph representation of Zachary's karate club network, where each vertex corresponds to one out of the 34 members of the club, and each edge indicates interaction between two members.  Note: this does not account for relative frequencies of interactions.\cite{Zach}}
\end{figure}

While predefined and comprehensively observed networks can clearly be visualized graphically and have tractable properties, problems arise when we begin to work with networks that are less tractable---where the structure of the network graph has properties which are not, or are only partially, observed.  This is where our statistical methods get interesting, as well-defined structure lends itself towards fairly straightforward analysis.  The solution to this is to construct network graphs in order to analyze observed structure and also allow for estimating structure which is not well-observed.\cite{Eric}

\subsection{Graphs and Linear Algebra}
Here, we present an intersection of graph theory and matrix algebra, predicated on an introductory understanding of each topic.  First, we define notation.  We will be using definitions consistent with B\'ona's.\cite{Walk}

\begin{itemize}

\item We define a graph structure $G(V,E)$ where $V$ is a set of vertices (nodes) and $E$ is a set of edges (links) such that for $u,v \in V$ and $u \neq v$ we have that $\{u,v\} \in E$ is an edge of $G$.

\item The \emph{order} of $G$ is denoted $N_v = \left|V\right|$, and the \emph{size} of $G$ is denoted $N_e = \left|E\right|$

\item We define \emph{adjacency} as follows.  For two vertices $u,v \in V$, we say $u$ is adjacent to $v$ if $\{u,v\} \in E$.  Similarly for edges, given $e_1, e_2 \in E$, we say $e_1$ is adjacent to $e_2$ if share an endpoint $w \in  V$.

\item We call a vertex $w \in V$ \emph{incident} to an edge $e_1 \in E$ if $w$ is an endpoint of $e_1$.

\end{itemize}
\begin{figure}
\centering
\caption{This is the adjacency matrix of Zachary's karate club network generated using MATLAB, where white blocks represent a $0$ in the $i,j$-th entry of {\bf{A}} and a black block represents a $1$ in the $i,j$-th entry of {\bf{A}}.\cite{Zach}}
\includegraphics[width=.5\textwidth]{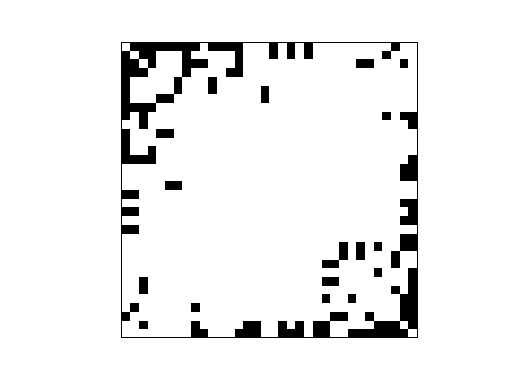}
\end{figure}

Having defined the requisite notation for the graph component of this section, we will examine the intersection now with matrix algebra.  Given $G(V,E)$ we define the adjacency matrix of {\bf{A}} of $G$ by $$A_i^j = \mathbbm{1}(i\sim_G j),$$ where $\mathbbm{1}(i\sim_G j)$ denotes the indicator function of the adjacency of vertices $i,j \in V$.  Of course, {\bf{A}} is $N_v$ x $N_v$, symmetric, and binary.  We can define the incidence matrix {\bf{B}} of $G$ similarly with the indicator function of the i-th vertex being an endpoint of the j-th row. $$B_i^j = \mathbbm{1}(i\sim_G j),$$ where {\bf{B}} is $N_v$ x $N_e$ and binary.

\section{Random Graph Models}
Having built a requisite background in network, graph, and matrix theory we will begin discussion around RGMs.  We will give a background on RGMs and give an examination of three different models; the Erd\H{o}s-R\'enyi, Watts-Strogatz, and Barab\'asi-Albert models.  

\subsection{Background of Random Graph Models}
Here, we establish a mathematical basis by which to proceed with our examination of RGMs.  Throughout this section we will be using notation and structure consistent with Kolaczyk.\cite{Eric}

Take for instance the problem of determining the significance of some structural characteristic, $\eta(\cdot)$, of a graph of observations, $G^{obs}$.  We seek to determine the significance of $\eta(G^{obs})$.  In accordance with Kolaczyk's method,\cite{Eric} we form a collection, $\mathcal{G}$, of random graphs.  We then compare $\{\eta(G) : G \in \mathcal{G}\}$ and our observed characteristic $\eta(G^{obs})$ and draw a conclusion based on the ``extreme-ness" of our observed characteristic in relation to our sampling of random graphs.  More formally\footnote{This definition, while dense, can be thought of as a the sum of a binary operator over the order of our graph family.  Namely, if we sum the binary operator which is 1 if our observed characteristic is less than some level $t$ for some graph $G_i \in \mathcal{G}$ and 0 else and then divide by the order of our graph family,, we then have the probability for our characteristic $\eta(G)$.} $$\mathbbm{P}_{\eta, \mathcal{G}}(t) = \frac{\sum\limits_{i=1}^{n} \mathbbm{1}(G_i \in \mathcal{G} : \eta(G_i) \leq t)}{\left|\mathcal{G}\right|}$$
Under this distribution, we note that unlikely $\eta(G^{obs})$ is evidence against a uniform draw from $\mathcal{G}$.\cite{Eric}

\subsection{The Erd\H{o}s-R\'enyi Model}
The classical random graph model is that originating out of several seminal papers from Paul Erd\H{o}s and Alfr\'ed R\'enyi.  Their model gives a collection of graphs $\mathcal{G}_{N_v, N_e}$ and then assigns some probability $\mathbbm{P}(\mathcal{G}_i) = {N \choose N_e}^{-1}$  for each $G_i \in \mathcal{G}$ where $N$ is the number of distinct vertex pairs.  

Further, we see that for some collection $\mathcal{G}_{N_v, p}$ where we define $p \in (0,1)$, we have $N_e\sim pN_v$ is equivalent to the above for $N_v \gg 0$.\cite{Gilbert}  These results were produced --- contemporaneously to Erd\H{o}s and R\'enyi's works --- by E.N. Gilbert.\cite{Gilbert}

To examine some properties of the Erd\H{o}s-R\'enyi model, we will look at the results produced pertaining to connectivity, degree distribution, and clustering.  Kolaczyk finds the level of connectivity of $G \in \mathcal{G}_{N_v, p}$ by relating $N_v$ and $p$.  By letting $p  = c/N_v$ for $c>0$, this gives us an expected density $p \sim \frac{1}{N_v}$.  This is indicative that $G$ is likely sparse for $N_v \gg 0$.

Erd\H{o}s and R\'enyi found the degree distribution generated by their model to follow a Poisson distribution as follows.

Define $\nu_r (G)$ as the number of vertices in graph $G$ with degree $r$.  Further, define $\Gamma_{n,N}$ as a  random graph with $N$ edges chosen from $n$ vertices.  From this, Erd\H{o}s and R\'enyi proved that, for $k \in \mathbb{Z}_{\geq 0}$ $$\lim_{n \to \infty}\mathbbm{P}(\nu_r(\Gamma_{n, N(n)})=k) = \frac{\lambda^{k}e^{-\lambda}}{k!}$$  Which we immediately recognize as a the probability mass function of a Poisson random variable with mean $\lambda$.\cite{ER}

From the properties of degree distribution and connectivity, we see that the Erd\H{o}s - R\'enyi model does a strong job of modeling the small-world property of networks, which is the average shortest path length, but falls short in modeling large-scale real-world networks.\cite{Eric}  As will be clear from the examination of clustering in the next section, this model also does not incorporate any sort of strong clustering properties.  In fact, for this model the clustering coefficient is $p$ by construction, and by definition $p$ approaches zero as $N_v$ tends towards infinity.

\subsection{Watts-Strogatz Model}
The shortcomings of the Erd\H{o}s-R\'enyi Model in modeling connectivity and clustering which mirrors that observed in large-scale real-world networks led to groundbreaking work by Duncan Watts and Steven Strogatz.  Watts and Strogatz recognized that Bernoulli random graphs --- that the Erd\H{o}s-R\'enyi model uses --- fall short in their lack of clustering, and so constructed a new random graph model which had the benefit of the small-world properties of the E-R graphs, but also incorporated strong clustering properties in their topology.

At this point we define $\lambda_G(v)$ as the number of triangles of $G$ --- i.e. the number of subgraphs with 3 edges and 3 vertices --- for which $v$ is a vertex.  Further, $\tau_G(v)$ is the number of subgraphs of $G$ with two edges and three vertices with $v$ as a vertex which is incident to both edges.\footnote{This is also equivalent to the number of 3-trees with $v$ as the root}  The clustering coefficient of a graph $G$ can then be defined as $$cl_T(G) = \frac{\lambda_G(v)}{\tau_G(v)}.$$

\begin{figure}
\centering
\caption{This is an example of the switching algorithm provided by Watts and Strogatz.  Each figure is based on a regular ring lattice, as seen on the left, where the edges are switched with probability zero.  While on the right we see a highly chaotic graph that does not appear to have nice properties.  In the intermediary stage, however, we note that the graph displays the small-world clustering with small shortest average path length.\cite{WS}}
\includegraphics[width=.75\textwidth]{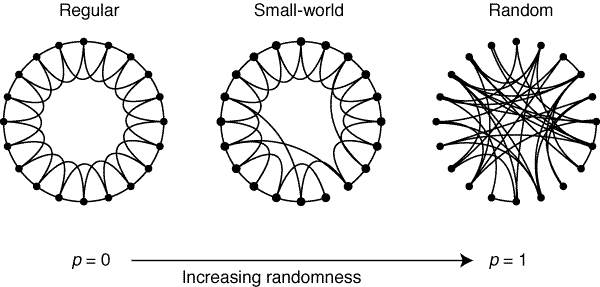}
\end{figure}

The Watts-Strogatz model is an algorithm which takes a regular ring lattice, as seen in Figure 3, and switch an edge between vertices with probability $p$.  The Watts-Strogatz results quantified characteristic path length $L(p)$ and clustering coefficient $C(p)$, and found that given, for $n$ vertices of degree $k$ with $n\gg k \gg \ln(n) \gg1$ and $k \gg \ln(n)$, the random graph generated will be connected.  In this model they found $L\sim n/2k \gg 1$ and $C\sim 3/4$ as $p \rightarrow 0$, and $L \approx L_{random} \sim \ln(n)/\ln(k)$ and $C \approx C_{random} \sim k/n \ll 1$ as $p \rightarrow 1$.  While this might appear to indicate an inverse relationship between $L$ and $C$, Watts and Strogatz contended that due to a handful of long-range links, or shortcuts, there is a wide range of $p$ for which $L(p)$ approximates $L_{random}$ with $C(p) \gg C_{random}$.\cite{WS}

\begin{figure}
\centering
\includegraphics[width=.7\textwidth]{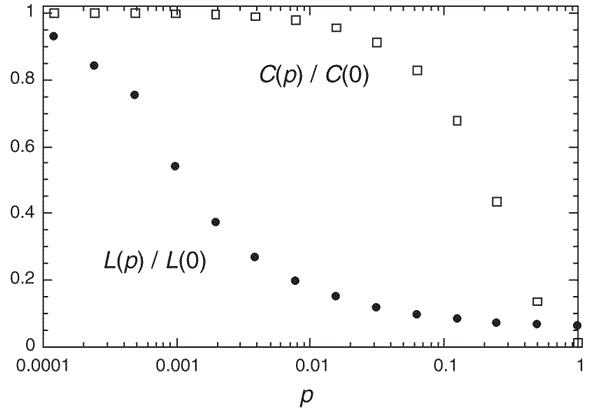}
\caption{Here we see that for a range of $p$, $C(p)$ is quite high while $L(p)$ is quite small.  Logarithmic horizontal scale was used due to the decay rate of $L(p)$ with the onset of the `small-world' effect, but locally as that effect takes place $C(p)$ remains fairly constant.\cite{WS}}
\end{figure}

From this work, which has been highly influential across disciplines, we see a way to reconcile the strengths and shortcomings of the classical random graph model.  We are left to note that, though having local clustering and small-world phenomena built into the topology of an RGM is remarkable, these models do not seem to account for the way in which real world connections are made nor the generation of new connections from addition of new nodes to the network.  Intuitively, new friendships are not evenly distributed across all people with some probability.  Rather, those with large social networks to begin with are those more likely to gain new friends---their corresponding vertex is more likely have an edge with a newly added node.  This is exactly the notion of a power-law distribution and preferential attachment.

\subsection{Barab\'asi-Albert Model}
Albert-L\'aszl\'o Barab\'asi and R\'eka Albert offered a model by which the topology of real-world networks could be built into random graphs.  Namely, real-world networks display a power-law distribution.  Stated explicitly, this is to say that the probability $P(k)$ that a vertex interacts with $k$ other vertices decays in the form $P(k) \sim k^{-\gamma}$.\cite{BA}  Using a few readily available datasets with well-known topology, this is evidenced.

\begin{figure}
\centering
\caption{Here, we have for {\bf A} the case of collaboration of actors in movies.  The probability of an actor having $k$ links is given by $P(k)$ which follows the power-law distribution for $\gamma \approx 2.3$.  Similarly, connections in the World Wide Web ({\bf B}) and the US electrical power grid ({\bf C}) behave with a power-law distribution.\cite{BA}}
\includegraphics[width=.75\textwidth]{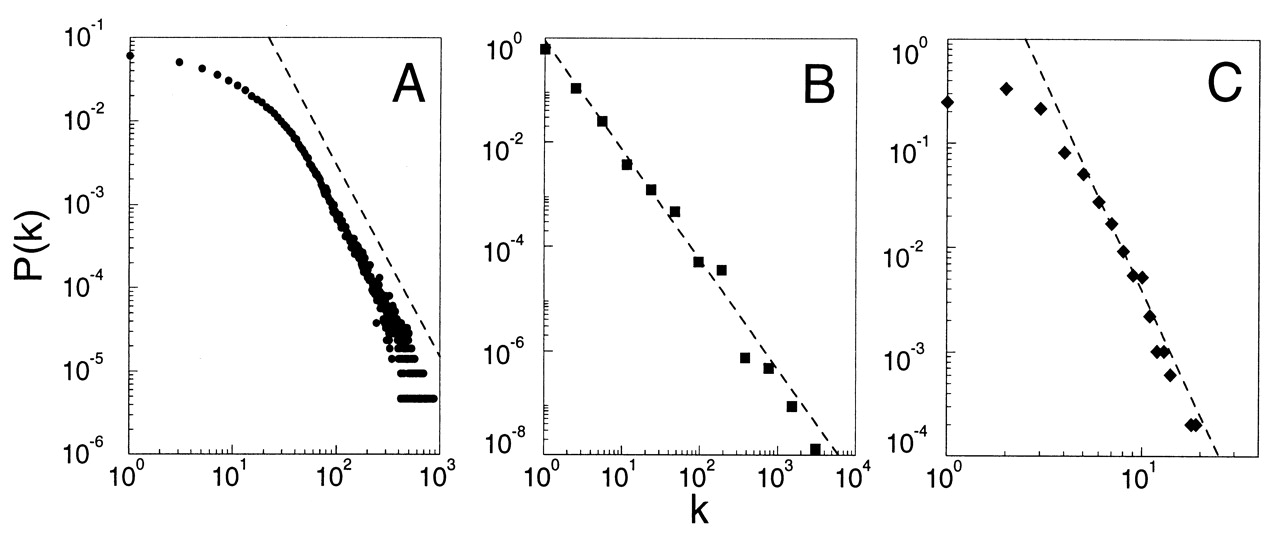}

\end{figure}

Barab\'asi and Albert built this scale-free distribution into their RGM by incorporating both growth in the network and preferential attachment.  Using the case of the actors, as shown in Figure 5, it is intuitive that a new actor is more likely to have a role with an experienced and established actor---one with a larger network.  By assigning a probability $\Pi$ to a new vertex's connection to an existing vertix $i$ which depends on the connectivity $k_i$ of that vertex, we have $\Pi(k_i) = k_i / \sum_j k_j$.  We take steps of time, $t$, and end with $t + m_0$ vertices and $mt$ edges.  This leads to the scale-invariant distribution that Barab\'asi and Albert sought to construct, using a model of growth and preferential attachment.  The rate of connectivity we get from this model is $\partial k/\partial t = k_i /2t$, and so our $t$-th moment is given by $k_i(t) = m(t/t_i)^{0.5}$ which leads to a probability density $P(k) = \partial P[k_i(t) \l k]/\partial k$ for a large interval of $t$ yields $$P(k) = \frac{2m^2}{k^3}.$$  This density yields $\gamma = 3$, which is independent of $m$---scale-invariance.

This gives us that which we sought --- to have a random graph model which incorporates growth and preferential connectivity in order to replicate the scale-invariance of real-world network growth models.

\section{Random Graph Matching and the Graph Matching Problem}
In this section we seek to place RGMs in the broader context of the literature, and will discuss some recent results pertaining to the Graph Matching Problem (GMP).  The focus of this section will be on correlated Bernoulli random graphs,  with a focus on the Erd\H{o}s-R\'enyi random graph family.  In particular, this will be focused on a few recent collaborative works by Vince Lyzinski (UMass Amherst) and a number of other researchers.  
\subsection{Alignment in Correlated Bernoulli Graphs}
 Recall our definition of the GMP in the introduction as finding structure-preserving alignment between the vertex sets of two graphs.  Lyzinski explores this in the case of the correlated Erd\H{o}s-R\'enyi random graph---two graphs with the same set of vertices, each generated in the form of Bernoulli($p$).  Before exploring the results, we will give some definitions.
\begin{flushleft}
For $\Phi : V(G_1) \rightarrow V(G_2)$ as the latent alignment function of the vertex sets of our correlated Erd\H{o}s-R\'enyi random graph we have\footnote{These definitions are directly from the paper ``Seeded Graph Matching for Correlated Erd\H{o}s-R\'enyi Graphs".\cite{Vince}}
\end{flushleft}
\begin{itemize}

\item Vertex $v \in V(G_1)$ is mismatched by graph matching if there exists a solution function $\psi$ such that $\Phi(v) \neq \psi(v)$.

\item The GMP provides a consistent estimate of $\Phi$ if the number of mismatched vertices goes to zero in probability as $\left|V(G_1)\right|$ tends to infinity.

\item The GMP provides a strongly consistent estimate of $\Phi$ if the number of mismatched vertices converges almost surely to zero as $\left|V(G_1)\right|$ tends to infinity.

\end{itemize}

For consistency in notation, we will use some slightly different notation than used in the paper we are referencing.  We construct an Erd\H{o}s-R\'enyi random graph using parameters $n = N_{v_i}$, $p \in (0,1)$ and $\rho \in [0,1]$ where $N_{v_1} = N_{v_2}$ for respective graphs $G_1$ and $G_2$. For vertex pairs $\{v,v'\}$, we let $\mathbbm{1}(\{v,v'\} \in E_i)$ be the random variable for the event $\{v,v'\} \in E_i$, which is i.i.d. Bernoulli($p$) with correlation coefficient $\rho$.   This formulation will be used in the context of the graph matching problem, as follows.

We define $\Pi$ to be the set of bijections from $V(G_1) \rightarrow V(G_2)$.  Then the disagreements under $\psi \in \Pi$ are represented as exactly $$\Delta(\psi) := \sum\limits_{\{v,v'\}\in {n\choose 2}} \mathbbm{1}\left(\mathbbm{1}\left(v \sim_{G_1}v'\right) \neq \mathbbm{1}\left(\psi(v) \sim_{G_2} \psi\left(v'\right)\right)\right)$$
The graph matching problem is exactly the set of bijections in $\Pi$ that minimizes edge disagreements, denoted $\Psi := \argmin_{\psi \in \Pi} \Delta(\psi).$  We can now present a component of the primary theoretical result from Lyzinski, Fishkind and Priebe's work.

\begin{theorem}
Fix a real number $\xi_1 < 1$ such that $p \leq \xi_1$.  Then for fixed $c_1, c_2, c_3, c_4 \in \mathbb{R}$, dependent only on $\xi_1$, we have:
\begin{enumerate}[i)]
\item If $\rho \geq c_1 \sqrt{\frac{\log(n)}{np}}$ and $p \geq c_2\frac{\log(n)}{n}$ then almost always $\Psi = \{\Phi\}$.
\end{enumerate}
\end{theorem}

This result states that given some relatively lenient conditioning on the correlation between $G_1$ and $G_2$, the graph matching problem provides a strongly consistent estimator of the latent alignment function for $G_1$ and $G_2$ which holds in both sparse and dense graph schemes.\cite{Vince}\footnote{This is directly related to the range of $p$, as discussed in our examination of the clustering coefficient}

Lyzinski et al. provide another result toward the graph matching problem in a paper that shows that alignment strength and total correlation are asymptotically equivalent in the case of Bernoulli correlated random graphs.  We define our alignment strength as a function of our above definition of disagreements between $G_1$ and $G_2$ as follows
$${\bf str}(G_1,G_2, \psi) := 1 - \frac{\Delta(\psi)}{\frac{1}{n!}\sum_{\psi' \in \Pi_n}\Delta(\psi')}$$
In the case of a known alignment, ${\bf str}(G_1,G_2, \psi)$ gives us a quantifier of the structural similarities of $G_1$ and $G_2$, while in the case of an unknown alignment we use the graph matching solution $\psi_{GM} \in \Phi$,\footnote{As defined in the above section} such that ${\bf str}(G_1,G_2, \psi_{GM})$ gives a quantifier of the structural similarity between $G_1$ and $G_2$.\cite{Vince2}

Lyzinski et al. then offer a novel function, $\rho_T$, which is asymptotically equivalent to ${\bf str}(G_1,G_2, \psi)$ and is a function of the intergraph correlation coefficient $\rho_e$ and a newly defined intragraph measure they call the heterogeneity correlation, denoted $\rho_h$.  Before giving the result, we must state a few more definitions, as given by Lyzinski et al.

Referring to the Bernoulli parameters ${p_{\{i,j\} \in {[n] \choose 2}}}$, let their mean be
\begin{center}
$\mu := \frac{1}{{n \choose 2}}\sum_{\{i,j\} \in {[n] \choose 2}} p_{i,j}$
\end{center}
and let their variance be
\begin{center}
$\sigma^2 :=\frac{1}{{n \choose 2}}\sum_{\{i,j\} \in {[n] \choose 2}} (p_{i,j} - \mu)^2.$
\end{center}
The heterogeneity coefficient is then defined to be 
\begin{center}
$\rho_h :=\frac{\sigma^2}{\mu(1-\mu)}.$
\end{center}
The total correlation, $\rho_T$, is then defined such that
\begin{center}
$1-\rho_T := (1-\rho_e)(1-\rho_h).$
\end{center}
The main theoretical result of Lyzinski et al. (2018) is exactly that $${\bf str}(G_1,G_2, I) - \rho_T \xrightarrow{a.s} 0$$ where $I$ is the identity bijection.\cite{Vince2}  This result opens a rich new avenue for approaching problems of matchability in the universal case, which are at present far from tractable.

\section{Discussion}
This paper is meant to serve as a reference on which to base further analysis of RGMs and to build an appreciation of some of the purposes they serve.  This is by no means comprehensive --- further reading on network analysis and random graph models is encouraged.  Exciting research is being done in this area, having found its beginnings with Erd\H{o}s and R\'enyi and continuing to recent results.  Worth noting, too, is that any universal solution which predicts matching will include the novel parameter $\rho_T$ which we saw from Lyzinski et al.\cite{Vince2}  Random graph models and graph matching prove to be an interesting intersection of matrix algebra, graph theory, and statistics.  If the reader has a preference to the matrix algebra representation, see \cite{Vince}.  While the GMP is far from tractable at present, we may soon see great breakthroughs utilizing the properties of the models that we discussed here.

\end{document}